\documentclass{gtart_h}

\def\ifplaintex{\expandafter\ifx\csname documentclass\endcsname\relax}

\def\gtp{{\mathsurround=0pt\it $\cal G\mskip-2mu$eometry \&\ 
$\cal T\!\!$opology $\cal P\!$ublications}}  

\def\Addressesr{\bigskip
{\small \parskip 0pt \leftskip 0pt \rightskip 0pt plus 1fil \def\\{\par}
\sl\theaddress\par
\medskip
\rm Email:\stdspace\tt\theemail\hfill\rm Received:\qua\receiveddate \par}}

\def\recd{{\small Received:\qua\receiveddate\ifx\reviseddate\relax
\else\qquad Revised:\qua\reviseddate\fi\par}} 

\def\lognumber#1{\def\thelognumber{#1}}
\def\volumenumber#1{\def\thevolumenumber{#1}}
\def\volumeyear#1{\def\thevolumeyear{#1}}
\def\papernumber#1{\def\thepapernumber{#1}}
\def\pagenumbers#1#2{\def\startpage{#1}\def\finishpage{#2}}
\def\published#1{\def\publishdate{#1}}

\def\received#1{\def\receiveddate{#1}}

\def\accepted#1{\def\accepteddate{#1}}

\def\asciiaddress#1{\def\theasciiaddress{#1}}

\long\def\asciiabstract#1{\long\def\theasciiabstract{#1}}

\let\\\par\let\thelognumber\relax\let\thevolumenumber\relax
\let\thepapernumber\relax\let\thevolumeyear\relax\let\startpage\relax
\let\finishpage\relax\let\publishdate\relax\let\receiveddate\relax
\let\reviseddate\relax\let\accepteddate\relax\let\theasciititle\relax
\let\theasciiauthors\relax\let\theasciiaddress\relax
\let\theasciiabstract\relax

\let\theasciiemail\relax

\ifplaintex
\font\logobig=cmssbx10 scaled 3836
\font\logomed=cmssbx10 scaled 2557
\else
\font\logobig=cmssbx10 scaled 4200
\font\logomed=cmssbx10 scaled 2800
\fi

\long\def\makeagttitle{   
\count0=\startpage
\agt\hfill      
\hbox to 45truept{\vbox to 0pt{\vglue -13truept{\logomed A\kern -.37em{\logobig 
T}\kern -.38em G}\vss}\hss}
\break
{\small Volume \thevolumenumber\ (\thevolumeyear)
\startpage--\finishpage\nl
Published: \publishdate}

\vglue .25truein

{\parskip=0pt\leftskip 0pt plus
1fil\def\\{\par\smallskip}{\Large\bf\thetitle}\par\medskip} \vglue
0.05truein

{\parskip=0pt\leftskip 0pt plus 1fil\def\\{\par}{\sc\theauthors}
\par\medskip}%
 
\vglue 0.03truein 

{\small\leftskip 25truept\rightskip 25truept{\bf Abstract}\stdspace\theabstract

{\bf AMS Classification}\stdspace\theprimaryclass
\ifx\thesecondaryclass\relax\else; \thesecondaryclass\fi\par
{\bf Keywords}\stdspace \thekeywords\par}\vglue 7truept

}   

\ifplaintex
\font\phead=cmsl9 scaled 950
\font\pnum=cmbx10 scaled 913
\font\pfoot=cmsl9 scaled 950
\headline{\vbox to 0pt{\vskip -4.5mm\line{\small\phead\ifnum
\count0=\startpage ISSN 1472-2739 (on-line) 1472-2747 (printed)
\hfill {\pnum\folio}\else\ifodd\count0\def\\{ }%
\ifx\theshorttitle\relax\thetitle\else\theshorttitle\fi\hfill{\pnum\folio}
\else\def\\{ and }{\pnum\folio}\hfill\ifx\theshortauthors\relax\theauthors
\else\theshortauthors\fi\fi\fi}\vss}}
\footline{\vbox to 0pt{\vglue 0mm\line{\small\pfoot\ifnum\count0=\startpage
\copyright\ \gtp\hfill\else
\agt, Volume \thevolumenumber\ (\thevolumeyear)\hfill\fi}\vss}}
\else
\font=cmsl9 scaled 1050
\font\lnum=cmbx10 
\font=cmsl9 scaled 1050
\makeatletter
\def\@oddhead{{\small\lhead\ifnum\count0=\startpage ISSN 1472-2739 
(on-line) 1472-2747 (printed)\hfill {\lnum\number\count0}\else\ifodd\count0
\def\\{ }\ifx\theshorttitle\relax \thetitle \else\theshorttitle\fi\hfill
{\lnum\number\count0}\else\def\\{ and }{\lnum\number\count0}
\hfill\ifx\theshortauthors\relax 
\theauthors\else\theshortauthors\fi\fi\fi}}\def\@evenhead{\@oddhead}
\def\@oddfoot{\small\lfoot\ifnum\count0=\startpage\copyright\ \gtp\hfill\else
\agt, Volume \thevolumenumber\ (\thevolumeyear)\hfill\fi}
\def\@evenfoot{\@oddfoot}
\makeatother
\fi
\let\maketitlepage\makeagttitle
\let\makeshorttitle\maketitlepage
\let\maketitle\maketitlepage

\newwrite\gtoutfile
\long\gdef\makeheadfile{  
{\def\\{, }\def\s{ }
\immediate\openout\gtoutfile head.xxx
\immediate\write\gtoutfile{Proxy-for: \ifx\theasciiauthors\relax
\theauthors\else\theasciiauthors\fi\s<\ifx\theasciiemail\relax\theemail\else\theasciiemail\fi>}
\immediate\write\gtoutfile{\noexpand\\}
\immediate\write\gtoutfile{Authors: \ifx\theasciiauthors\relax
\theauthors\else\theasciiauthors\fi}
{\def\\{ }\immediate\write\gtoutfile{Title: \ifx\theasciititle\relax
\thetitle\else\theasciititle\fi}}
\immediate\write\gtoutfile{Subj-class: GT or SG, GR etc}
\immediate\write\gtoutfile{MSC-class: \theprimaryclass\ifx\thesecondaryclass\relax\else, \thesecondaryclass\fi}
\immediate\write\gtoutfile{Journal-ref: Algebr. Geom. Topol. \thevolumenumber\s
(\thevolumeyear) \startpage-\finishpage}
\immediate\write\gtoutfile{Comments: Published by Algebraic and
Geometric Topology at}
\immediate\write\gtoutfile{\s\s\s  http://www.maths.warwick.ac.uk/agt/AGTVol\thevolumenumber/agt-\thevolumenumber-\thepapernumber.abs.html}
\immediate\write\gtoutfile{\noexpand\\}
\immediate\write\gtoutfile{}
\ifx\theasciiabstract\relax
\immediate\write\gtoutfile{\theabstract}\else
\immediate\write\gtoutfile{\theasciiabstract}\fi
\immediate\write\gtoutfile{}
\immediate\write\gtoutfile{\noexpand\\}
\immediate\write\gtoutfile{}
\immediate\closeout\gtoutfile}}  

\def\maketitlepage{\makeagttitle\makeheadfile}
\let\makeshorttitle\maketitlepage
\let\maketitle\maketitlepage

\lognumber{18}
\volumenumber{4}
\volumeyear{2004}
\papernumber{18}
\published{29 May 2004}
\pagenumbers{333}{346}
\received{24 November 2003}
\accepted{29 May 2004}

\def\sh#1{\goodbreak{\large\bf #1}\nobreak
\addcontentsline{toc}{subsection}{#1}}

\usepackage{amsmath,amscd,amssymb,enumerate}

\swapnumbers
\theoremstyle{plain}
\newtheorem{theorem}{Theorem}[section]
\newtheorem{lemma}[theorem]{Lemma}
\newtheorem{corollary}[theorem]{Corollary}
\newtheorem{proposition}[theorem]{Proposition}

\theoremstyle{definition}
\newtheorem{definition}[theorem]{Definition}

\newtheorem{remark}[theorem]{Remark}

\newcommand{\reals}{\mathbb{R}}
\newcommand{\complexs}{\mathbb{C}}

\newcommand{\integers}{\mathbb{Z}}
\newcommand{\rationals}{\mathbb{Q}}

\DeclareMathOperator{\id}{id}

\newcommand{\boundedops}{\mathcal{B}}

\newcommand{\tensor}{\otimes}
\newcommand{\into}{\hookrightarrow}

\newcommand{\iso}{\cong}

\DeclareMathOperator{\im}{im}      

\DeclareMathOperator*{\dirlim}{dirlim}

\newcommand{\forget}[1]{}

\def  \nuint {\raise10pt\hbox{$\nu$}\kern-6pt\int}

\def \Sp {{\cal S}}

\newcommand\Di{D\kern-6pt/}
\newcommand\cDi{{\mathcal D}\kern-6pt/}
\newcommand\spi{S\kern-6pt/}
\newcommand \cspi{\Sp\kern-6pt/}

\def \s {\smallskip}

\def \cal {\mathcal}

{\catcode`@=11\global\let\c@equation=\c@theorem}

\renewcommand{\theenumi}{(\arabic{enumi})}

\allowdisplaybreaks[2]

\begin{document}

\title{Real versus complex K-theory using Kasparov's\\bivariant
  KK-theory}
\shorttitle{Real versus complex K-theory}
\author{Thomas Schick}

\address{Fachbereich Mathematik, Georg-August-Universit{\"a}t
G{\"o}ttingen, Germany}
\asciiaddress{Fachbereich Mathematik, Georg-August-Universitaet
Goettingen, Germany}

\email{schick@uni-math.gwdg.de}
\urladdr{http://www.uni-math.gwdg.de/schick}

\begin{abstract}
  In this paper, we use the KK-theory of Kasparov to prove exactness
  of sequences relating the K-theory of a real $C^*$-algebra and of
  its complexification (generalizing results of Boersema).

    We use this to relate the real version of the Baum-Connes conjecture
  for a discrete group to its complex counterpart. In particular, the
  complex Baum-Connes assembly map is an
  isomorphism if and only if the real one is, thus reproving a result
  of Baum and Karoubi. After inverting
  $2$, the same is true for
  the injectivity or surjectivity part alone.
\end{abstract}

\asciiabstract{%
In this paper, we use the KK-theory of Kasparov to prove exactness of
  sequences relating the K-theory of a real C^*-algebra and of its
  complexification (generalizing results of Boersema). We use this to
  relate the real version of the Baum-Connes conjecture for a discrete
  group to its complex counterpart. In particular, the complex
  Baum-Connes assembly map is an isomorphism if and only if the real
  one is, thus reproving a result of Baum and Karoubi. After inverting
  2, the same is true for the injectivity or surjectivity part alone.}



\maketitle

\section{Motivation}
\label{sec:real-versus-complex}

In the majority of available sources about the subject, complex
$C^*$-algebras and Banach algebras
and their K-theory is studied. However, for geometrical reasons, the
real versions also play a prominent role. 

Before describing the results of this paper, we want to give the
geometric motivation why both variants are necessary.

\begin{enumerate}
\item Real K-Theory (meaning K-theory of real $C^*$-algebra) is more powerful since it contains additional
  information. Most notably this can be seen at Hitchin's
  $\integers/2$-obstructions to positive scalar curvature in
  dimensions $8k+1$ and $8k+2$ \cite{MR50:11332}. They take values in
  $KO_{j}(\reals)$ for $j=1,2$. Related to this is
  the fact that there are $8$ different
  groups, and not just $2$, since real K-theory does not have the
  $2$-periodicity of complex K-theory, but is $8$-periodic.

  In particular, we mention the following result of Stephan Stolz: if
  the real Baum-Connes map $\mu_{\reals,red}\colon
  RKO^\Gamma_*(\underline{E}\Gamma)\to KO_*(C^*_{\reals,red}\Gamma)$
  is injective, then the stable Gromov-Lawson-Rosenberg conjecture is
  true for $\Gamma$. This means that a spin manifold with fundamental
  group $\Gamma$ stably admits a metric with positive scalar curvature
  if and only if the Mishchenko-Fomenko index of its Dirac operator vanishes.
\item Unfortunately, a real structure of some kind is needed to define indices in real
  K-theory. In particular, there is no good way to define a (higher) real
  index of the signature operator in dimension $4k+2$.
\end{enumerate}

This explains why for the Dirac operator, and therefore for the study
of metrics of positive scalar curvature on spin manifolds, one traditionally uses real
K-theory, whereas complex K-theory is used for the signature operator
and the study of higher signatures.

This issue came up in the paper \cite{piazza:_bordis} of Paolo Piazza
and the author, where we studied both the signature operator and the
spin Dirac operator.

\section{Real versus complex K-theory}
\label{sec:real-versus-complex-1}

In this paper, we give a
theoretical comparison of real and complex K-theory. The results of this short note are
essentially ``folklore'' knowledge. Early results  date back to
\cite{A1}. However, there only the special case of commutative
$C^*$-algebras (in other words, spaces) is considered.

General results about the relation between real and complex K-theory
are proved by Max Karoubi in \cite{Karoubi2}, using some modern homotopy theory. The
results of \cite{Karoubi2} are applied in \cite{BKR} by Paul Baum and
Max Karoubi to prove that,
for discrete groups, the
complex Baum-Connes conjecture implies the real Baum-Connes
conjecture. Their proof is based, apart from \cite{Karoubi2}, on the
interpretation of the Baum-Connes map as a connecting homomorphism as
explained by Roe in \cite{MR97h:58155}. Our results are related to and in part
equal to their results. We use, however, a different method entirely
embedded in (real) bivariant K-theory (i.e.~KK-theory), as developed
by Kasparov (compare e.g.~\cite{MR88j:58123} and \cite{MR81m:58075}).

In the non-equivariant setting, the exact sequences stated below relating real and
complex KK-theory are established (with similar methods) by Boersema
in \cite{math.OA/0302335} and \cite{MR1935138}. His united KK-theory
can be extended to the equivariant setting and then the framework of
the so called acyclic CRT-modules and there properties as introduced
by Bousfield \cite{MR92d:55003} could be used to give another proof of
the equivalence of the real and complex Baum-Connes conjecture.

To keep the paper self contained, we reprove a number of results which
(modulo extension to the equivariant case) can be found in
\cite{math.OA/0310209,math.OA/0302335}. 

We prove the following theorems.

\begin{theorem}\label{theo:long_exact_K-theory}
  Let $A$ be a separable \emph{real} $\sigma$-unital $C^*$-algebra and
  $A_\complexs:=A\tensor \complexs$. Then there is a long exact sequence in K-theory of
  $C^*$-algebras
  \begin{equation}\label{eq:long_exact_K_sequence}
   \cdots \to KO_{q-1}(A)\xrightarrow{\chi} KO_{q}(A)\xrightarrow{c}
    K_{q}(A_\complexs) \xrightarrow{\delta} KO_{q-2}(A)\to\cdots
  \end{equation}
  Here, $c$ is complexification, $\chi$ is multiplication by the
  generator $\eta\in KO_1(\reals)\iso \integers/2$ (in particular
  $\chi^3=0$), and $\delta$ is the
  composition of the inverse of multiplication with the Bott element
  in $K_2(\complexs)$ with ``forgetting the complex structure''.
\end{theorem}

\begin{remark}
  Real and complex $C^*$-algebras and their K-theory are connected by
  ``complexification'' and ``forgetting the complex structure''. We
  use these terms throughout, precise definitions are given in
  Definitions \ref{def:complexification} and \ref{def:forgetting complex structure}.
\end{remark}

\begin{corollary}\label{corol:split_exact_K-theory}
  In the situation of Theorem \ref{theo:long_exact_K-theory},
  if we invert $2$, in particular if we tensor with $\rationals$, the
  sequence splits into short \emph{split} exact sequences
  \begin{equation}\label{eq:short_exact_K_sequence}
    0 \to 
    KO_{q}(A) \tensor\integers[\frac{1}{2}] \xrightarrow{c} 
    K_{q}(A_\complexs)\tensor \integers[\frac{1}{2}]
    \xrightarrow{\delta} KO_{q-2}(A)\tensor \integers[\frac{1}{2}]\to
    0 .
  \end{equation}
\end{corollary}
\begin{proof}
  We obtain short exact sequences because $2\eta=0$, i.e.~the image of $\chi$
  (and therefore the kernel of $c$) in
  \eqref{eq:long_exact_K_sequence} consists of $2$-torsion.

  The sequence is
 {split} exact, with split being given by ``forgetting the
  complex structure'' $K_*(A_\complexs)\to KO_*(A)$, since the
  composition of ``complexification'' with ``forgetting the complex
  structure'' induces multiplication with $2$ in $KO_*(A)$, i.e.~an
  automorphism after inverting $2$. For more
  details, compare Definition \ref{def:complexification}, Definition
  \ref{def:forgetting complex structure} and Lemma
  \ref{lem:complexification_and forgetting}.
\end{proof}

\begin{theorem}\label{theo:long_exact_K_homology}
 Assume that $\Gamma$ is a discrete group and $X$ is a proper
$\Gamma$-space. Let $B$ be a separable real $\sigma$-unital
$\Gamma$-$C^*$-algebra.  Then we have a long exact sequence in equivariant representable
K-homology with coefficients in $B$ (defined e.g.~via Kasparov's KK-theory)
  \begin{equation}\label{eq:long_exact_equiv_space_K_homology_sequence}
   \cdots\to  RKO^{\Gamma}_{q-1}(X;B)\xrightarrow{\chi} RKO^{\Gamma}_{q}(X;B)\xrightarrow{c}
    RK^{\Gamma}_{q}(X;B_\complexs) \xrightarrow{\delta} RKO^{\Gamma}_{q-2}(X;B)\cdots
  \end{equation}
  Here, $c$ is again complexification, and $\chi$ is given by
  multiplication with the generator in $KO_1(pt)=\integers/2$, i.e.~$\chi^3=0$. $\delta$
  is the composition of (the inverse of) the complex Bott periodicity
  isomorphism with ``forgetting the complex structure''. 
\end{theorem}

\begin{corollary}\label{corol:split_short_K-homology}
  In the situation of Theorem \ref{theo:long_exact_K_homology}, after
  inverting $2$, in particular after tensor product with
  $\rationals$, we obtain split short exact sequences
  \begin{equation}\label{eq:short_exact_equiv_K_homology_seq}
   0\to 
   KO^{\Gamma}_{q}(X;B) \tensor \integers[\frac{1}{2}]\xrightarrow{c} 
    K^{\Gamma}_{q}(X;B_\complexs) \tensor \integers[\frac{1}{2}]\xrightarrow{\delta}
    KO^{\Gamma}_{q-2}(X;B)\tensor  \integers[\frac{1}{2}] \to 0.
  \end{equation}
\end{corollary}
\begin{proof}
  Compare the proof of Corollary \ref{corol:split_exact_K-theory}.
\end{proof}

\begin{theorem}\label{theo:real_versus_complex_BC}
  Let $\Gamma$ be a discrete group. Consider the special case of
  Theorem \ref{theo:long_exact_K-theory} where
  $A=C^*_{\reals,red}(\Gamma;B)$ is the crossed product of $B$ by
  $\Gamma$, and the special case of Theorem
  \ref{theo:long_exact_K_homology} where $X=\underline{E}\Gamma$, the
  universal space for proper $\Gamma$-actions. We have
  (Baum-Connes) index maps
  \begin{align}\label{eq:index_maps}
    \mu_{red} &\colon RK^\Gamma_p(\underline{E}\Gamma;B_\complexs)\to
    K_p(C_{red}^*(\Gamma;B_\complexs));\\
    \mu_{\reals,red} &\colon
    RKO^\Gamma_p(\underline{E}\Gamma;B) \to K_p(C^*_{\reals,red}(\Gamma;B)).
  \end{align}
  Using the canonical identification
  $C^*_{\reals,red}(\Gamma;B)_{\complexs}=
  C^*_{red}(\Gamma;B_\complexs)$, the index 
  maps \eqref{eq:index_maps} commute with the maps in the long exact
  sequences \eqref{eq:long_exact_K_sequence} and
  \eqref{eq:short_exact_equiv_K_homology_seq}. 
\end{theorem}

\begin{corollary}\label{corol:real_versus_complex BC}
  The real Baum-Connes conjecture is true if and only if the complex
  Baum-Connes conjecture is true, i.e.~$\mu_{red}$ of
  \eqref{eq:index_maps} is an isomorphism if and only if
  $\mu_{\reals,red}$ is an isomorphism.

  After inverting $2$, in particular after tensoring with $2$,
  injectivity and surjectivity are separately equivalent in the real
  and complex case, i.e.~in
  \begin{align*}
    \mu_{red}\tensor \id_{\integers[\frac{1}{2}]} & \colon
      RK^\Gamma_p(\underline{E}\Gamma;B_\complexs)\tensor\integers[\frac{1}{2}]\to 
    K_p(C_{red}^*(\Gamma;B_\complexs))\tensor \integers[\frac{1}{2}];\\
 \mu_{\reals,red}\tensor\integers[\frac{1}{2}] &\colon
    RKO^\Gamma_p(\underline{E}\Gamma;B)\tensor\integers[\frac{1}{2}] \to
    K_p(C^*_{\reals,red}(\Gamma;B))\tensor\integers[\frac{1}{2}],
  \end{align*}
  one of the maps is injective for all $p$ if and only if the other maps is
  injective for all $p$, and is surjective for all $p$ if and only if
  the other map is surjective for all $p$.
\end{corollary}
\begin{proof}
  Using the long exact sequences \eqref{eq:long_exact_K_sequence} and
  \eqref{eq:long_exact_equiv_space_K_homology_sequence} and the
  $5$-lemma, if $\mu_{\reals,red}$ is an isomorphism then also
  $\mu_{red}$ is an isomorphism. For the converse, we use the
  algebraic Lemma \ref{lem:algebraic_iso} and the fact that
  $\chi^3=0$.

  After inverting $2$, the long exact sequences split into short exact
  sequences, and consequently we can deal with injectivity and
  surjectivity separately, using e.g.~the general form of the
  $5$-lemma \cite[Proposition 1.1]{CaEi}.
\end{proof}

\begin{theorem}\label{theo:variations}
  Corresponding results to the ones stated above hold if we replace
  the reduced $C^*$-algebras with the maximal ones (and the reduced
  index map with the maximal assembly map).

  Corresponding results also hold if we replace the classifying space
  for proper actions  $\underline{E}\Gamma$ with the classifying space
  for free actions $E\Gamma$.
If $B\Gamma:=E\Gamma/\Gamma$ is a finite CW-complex,
  and $B=\reals$, then $RKO_p(B\Gamma,B) = KO_p(B\Gamma)$ is the real
  K-homology of the space $B\Gamma$. 
We get the new index map as
  composition of the index map of Theorem
  \ref{theo:real_versus_complex_BC} with a canonical map
  $RKO^\Gamma_p(E\Gamma;B)\to RKO^\Gamma_p(\underline{E}\Gamma;B)$.
\end{theorem}

\begin{remark}
  Of course, in Theorem \ref{theo:variations}, the assembly map will in
  many cases \emph{not} be an isomorphism ---whereas no
  example is known such that the assembly map of Theorem
  \ref{theo:real_versus_complex_BC} is not an isomorphism. In Theorem
  \ref{theo:variations} we only claim that it is an isomorphism for
  the real version if and only if is an isomorphism for the complex version.
\end{remark}

The long exact sequences of Theorem \ref{theo:long_exact_K-theory} and
Theorem \ref{theo:long_exact_K_homology} are special cases of the
following bivariant theorem.

\begin{theorem}\label{theo:general_KK sequence}
  Let $\Gamma$ be a discrete group and $A$, $B$ separable real
  $\sigma$-unital $\Gamma$-$C^*$-algebras. Then there is a long exact sequence
  \begin{equation}\label{eq:bivariant_long_sequence}
   \cdots\! \to\! KKO^{\Gamma}_{q-1}(A;B)\xrightarrow{\chi} KKO^{\Gamma}_{q}(A;B)\xrightarrow{c}
    KK^{\Gamma}_{q}(A_\complexs;B_\complexs) \xrightarrow{\delta}
    KKO^{\Gamma}_{q-2}(A;B)\cdots 
  \end{equation}
  Here, $\chi$ is given by Kasparov product with the
  generator of $KKO^\Gamma_1(\reals,\reals)=\integers/2$, $c$ is given
  by \emph{complexification} as defined in Definition
  \ref{def:complexification}, and $\delta$ is the composition of the
  inverse of the complex Bott periodicity isomorphism with
  ``forgetting the complex structure'' as defined in Definition
  \ref{def:forgetting complex structure}.

  In particular, $2\eta=0$, $2\chi=0$, and $\chi^3=0$.
\end{theorem}


\section{Proofs of the theorems}
\label{sec:proofs-theorems}

Note first that Theorem \ref{theo:long_exact_K-theory} and Theorem
\ref{theo:long_exact_K_homology} indeed are special cases of
Theorem \ref{theo:general_KK sequence}. For Theorem
\ref{theo:long_exact_K-theory} we simply have to take $\Gamma=\{1\}$, $A=\reals$ (and
then $B$ of Theorem \ref{theo:general_KK sequence} is $A$ of Theorem
\ref{theo:long_exact_K-theory}). For
Theorem \ref{theo:long_exact_K_homology} let first $Y$ be a
$\Gamma$-compact $\Gamma$-invariant subspace of $X$, and set
$A=C_0(Y)$. By definition,
\begin{equation*}
  RKO_{p}^\Gamma(X;B) = \dirlim KKO_p^{\Gamma}(C_0(Y), B),
\end{equation*}
where the (direct) limit is taken over all $\Gamma$-compact subspaces
of $X$. The corresponding sequence for each $Y$ is exact. Since the
direct limit functor is exact, the same is true for the sequence
\eqref{eq:long_exact_equiv_space_K_homology_sequence}.

We therefore only have to prove Theorem \ref{theo:general_KK
  sequence}, Theorem \ref{theo:real_versus_complex_BC} and Lemma
\ref{lem:algebraic_iso} (which was used in the proof of Corollary
\ref{corol:real_versus_complex BC}).

\sh{An algebraic lemma}

\begin{lemma}\label{lem:algebraic_iso}
  Assume that one has a commutative diagram of abelian groups with
  exact rows which are $3$-periodic:
  \begin{equation}
    \label{eq:abstract_ex}
    \begin{CD}
      \cdots @>>> A @>{\chi}>> A @>{c}>> B @>{\delta}>> A @>>>
      \cdots\\
      && @VV{\mu_A}V @VV{\mu_A}V @VV{\mu_B}V @VV{\mu_A}V\\
      \cdots @>>> U @>{\chi_U}>> U @>{c_U}>> V @>{\delta_V}>> U @>>> \cdots
    \end{CD}
 \end{equation}
 Let $\chi$ and $\chi_U$ be endomorphisms of finite order. Then $\mu_B$ is an
 isomorphism if and only if the same is true for $\mu_A$.
\end{lemma}
\begin{proof}
  If $\mu_A$ is an isomorphism so is $\mu_B$ by the $5$-lemma.

  Perhaps the most elegant way to prove the converse is to observe that the
  rows from exact couples in the sense of \cite[Section
  2.2.3]{McCle}. Consequently, we get derived commutative diagrams of
  abelian groups with exact rows which are $3$-periodic:
  \begin{equation}
    \label{eq:abstract_ex_der}
    \begin{CD}
      \cdots @>>> \chi^n(A) @>{\chi|}>> \chi^n(A) @>{c_n}>> B_n
      @>{\delta_n}>> \chi^n(A) @>>>
      \cdots\\
      && @VV{\mu_A|}V @VV{\mu_A|}V @VV{(\mu_B)_n}V @VV{\mu_A|}V\\
      \cdots @>>> \chi_U^n(U) @>{\chi_U}|>> \chi_U^n(U) @>{(c_U)_n}>>
      V_n @>{(\delta_V)_n}>> \chi_U^n(U) @>>> \cdots
   \end{CD}
 \end{equation}
 Here, $\chi^n(A)$ is the image of $A$ under the $n$-fold iterated map
 $\chi$. One defines inductively $B_n:= \ker(c_{n-1}\circ
 \delta_{n-1})/\im(c_{n-1}\circ\delta_{n-1})$; this is a certain homology of
 $B_{n-1}$, $c_n$ is the composition of $c_{n-1}|_{\chi^n(A)}$ with the
 projection map, and $\delta_n$ and $(\mu_B)_n$ are the maps induced by
 $\delta_{n-1}$ and $(\mu_{B})_{n-1}$, respectively, which one proves are
 well defined on homology classes.

 In particular, $(\mu_B)_n$ is an isomorphism if $\mu_B$ is (an
 isomorphism of chain complexes induces an isomorphism on homology).

 We prove now by reverse induction that $\mu_A|\colon
 \chi^n(A)\to\chi_U^n(V)$ is an isomorphism for each $n$. Since $\chi$
 and $\chi_U$ have finite order, there images are eventually zero, so
 the assertion is true for $n$ large enough.

 Under the assumption that $\mu_A|_{\chi^n(A)}$ is an isomorphism, we have to
   prove the same for $\mu_A|\colon \chi^{n-1}(A)\to
   \chi_U^{n-1}(U)$. For this, consider the commutative diagram with
   exact rows
   \begin{equation}\label{eq:expalin-diag}
     \begin{CD}
        \chi^{n-1}(A)/\chi^n(A) @>{c_{n-1}}>> B_{n-1} @>{\delta_{n-1}}>> \chi^{n-1}(A) @>{\chi}>>
       \chi^n(A) @>>> 0\\
        @VV{\mu_A|}V @VV{(\mu_B)_{n-1}}V @VV{\mu_A|_{\chi^{n-1}(A)}}V
        @VV{\mu_A|_{\chi^n(A)}}V @VVV\\
       \chi_U^{n-1}(U)/\chi_U^n(U) @>{{c_U}_{n-1}}>> V_{n-1}
       @>{(\delta_V)_{n-1}}>> \chi_U^{n-1}(U) @>{{\chi_U}}>>
       \chi_U^n(U) @>>> 0\\
     \end{CD}
   \end{equation}
   obtained by cutting the long exact sequence
   \eqref{eq:abstract_ex_der}. We have just argued that $(\mu_B)_{n-1}$
   is an isomorphism since $\mu_B$ is one by assumption, and by the induction
   assumption $\mu_A|_{\chi^n(A)}$ is an
   isomorphism. By the $5$-lemma \cite[Proposition 1.1]{CaEi}
   $\mu_A|_{\chi^{n-1}(A)}$ therefore is onto. This implies that the
   leftmost vertical map in \eqref{eq:expalin-diag} also is onto. Now
   we can use the $5$-lemma \cite[Proposition 1.1]{CaEi} again to
   conclude that $\mu_A|_{\chi^{n-1}(A)}$ also is injective.

   Induction concludes the proof.
\end{proof}


\begin{remark}
  The proof of ``injectivity implies injectivity'' in Corollary
  \ref{theo:real_versus_complex_BC} follows from the fact that, after
  inverting $2$, 
  \begin{equation*}
K_p(B_\complexs)\tensor \integers[\frac{1}{2}]\iso
  (KO_p(B)\oplus KO_{p-2}(B))\tensor\integers[\frac{1}{2}]
\end{equation*}
for any real $C^*$-algebra $B$ in a
  natural way, with a similar assertion for the left hand side of the
  Baum-Connes assembly map.

  In particular, injectivity or surjectivity, respectively, in degree $p$ for the
  complex Baum-Connes map is (after inverting $2$) equivalent to
  injectivity or surjectivity, respectively, in the two degrees $p$
  and $p-2$ for the real Baum-Connes map.

  We should note that not only the proof of this assertion in Theorem
  \ref{theo:real_versus_complex_BC} does not work if $2$ is not
  inverted. Worse: the underlying algebraic statement is actually
  false. The easiest
  example is given by the short exact sequence for $K$-theory of a
  point. If we tensor this with $\integers[1/2]$, it remains
  exact. The natural map $M\to M\tensor \integers[1/2]$ connects the
  original exact sequence with the new exact sequence. For complex
  $K$-theory, the relevant maps are the inclusion $\integers\into
  \integers[1/2]$ and $0\to 0$. In particular, this is injective in
  all degrees. For real K-theory, we also get (in degrees $1$ and $2$
  mod $8$) the map $\integers/2\to 0$,
  i.e.~the map here is not injective in all degrees.
\end{remark}

\sh{Exterior product with ``small'' KK-elements}

\begin{definition}\label{def:small_homomos}
  In the sequel, we will frequently encounter homomorphisms $f\colon
  KKO^\Gamma_p\Gamma(A,B\tensor M_1)\to KKO^\Gamma_{p+l}(A,B\tensor
  M_2)$, where $M_1$, $M_2$ are ``elementary'' $C^*$-algebras (with
  trivial $\Gamma$-action),
  e.g.~$M_i\in\{\reals,\complexs,M_2(\reals),\cdots\}$. In most cases,
  $f$ will be induced by exterior Kasparov product with an element
  $[f]\in KK_l(M_1,M_2)$.
  
  Such a homomorphism will be called \emph{small}, or given by
  \emph{Kasparov product with a small element}. It is clear that the
  composition of small homomorphisms is again a small homomorphism.
\end{definition}

\sh{Complexification and forgetting the real structure}


\begin{definition}\label{def:complexification}
  Let $A$ be a real $\Gamma$-$C^*$-algebra with complexification
  $A_\complexs:=A\tensor \complexs$. Note that $A_\complexs$ can also be viewed as a real
  $\Gamma$-$C^*$-algebra, with a canonical natural inclusion $A\into
  A_\complexs$. This map, and the maps it induces on K-theory are
  called ``complexification'' and denoted by $c$.

 For $A$ and $B$
  separable real $\Gamma$-$C^*$-algebras, $A$ unital, the map induced
  by $c\colon A\into
  A_\complexs$ can be composed with the isomorphism of Proposition
  \ref{prop:compare_real_and_complex}, to get
  \begin{equation*}
  c_\complexs\colon KKO^\Gamma_n(B,A) \xrightarrow{[c]} KKO^\Gamma_n(B,A_\complexs) \iso
  KK^\Gamma_n(B_\complexs,A_\complexs).
\end{equation*}
This is what we call the complexification homomorphism in KK-theory,
it induces corresponding maps in K-theory and K-homology. Note that
$c$ is a small homomorphism, induced from $[c\colon
\reals\to\complexs] \in KKO_0(\reals,\complexs)$.
\end{definition}

\begin{definition}\label{def:forgetting complex structure}
  Let $A$ be a real separable $\Gamma$-$C^*$-algebra, $A_\complexs$ its
  complexification. We have a canonical natural inclusion
  $A_\complexs\into M_2(A)$, using the usual inclusion
  $i\colon\complexs\to M_2(\reals)$. If $A$ is $\sigma$-unital
  and $B$ is another separable $\Gamma$-$C^*$-algebra, using Morita
  equivalence, we get the induced homomorphisms in K-theory
  \begin{equation*}
   f_{\complexs}\colon KK^\Gamma_n(B_\complexs,A_\complexs)\iso KKO_n^\Gamma(B,
    A_\complexs)\xrightarrow{f}
    KKO_n^\Gamma(B,M_2(A))\xrightarrow[{M}]{\iso} KKO_n^\Gamma(B,A),
  \end{equation*}
  called ``forgetting the complex structure''. Note that $f$ is a
  small homomorphism, induced by $[i\colon\complexs\into M_2(\reals)]
  \in KKO_0(\complexs,M_2(\reals))$. Also $M$ is a small homomorphism,
  as explained in the proof of Lemma \ref{lem:complexification_and
    forgetting}. We define $f_\reals:=[i]\bullet [M_\reals]\in
  KKO_0(\complexs,\reals)$ to be the corresponding composition of
  small KKO-elements.
\end{definition}

\begin{lemma}\label{lem:complexification_and forgetting}
  In the situation of Definitions \ref{def:complexification} and
  \ref{def:forgetting complex structure}, the composition of first
  complexification and then forgetting the complex structure is
  multiplication by~$2$.
\end{lemma}
\begin{proof}
  By definition, this composition is the small homomorphism given by
  exterior Kasparov product with  $[i\colon \reals\into M_2(\reals)]$
  ($i$ the diagonal inclusion)
  composed with the small Morita equivalence homomorphism $[M_\reals] \in
  KKO_0(M_2(\reals),$ $\reals)$.
  It is known that $[i]\bullet [M_\reals] = 2\in
  KKO_0(\reals,\reals)$, which by associativity implies the claim. For
  a short KK-theoretic proof, observe that $[M_\reals] =
  [\reals^2\oplus 0, 0]\in KKO_0(M_2(\reals,\reals)$, with the
  obvious left $M_2(\reals)$ and right $\reals$-module structure on $\reals^2$, and
  with operator $0$. On the other hand, $[i]=(M_2(\reals)\oplus 0,
  0)$. Since both operators in our
  representatives are $0$ we get (compare e.g.~\cite{Bla}) 
  \begin{equation*}
    [i]\bullet [M_\reals] = [\reals^2\oplus 0, 0] = 2[\reals\oplus 0,0]\in KK^\Gamma_0(\reals,\reals).
  \end{equation*}
  Since $[\reals\oplus 0,0]=1\in KK^\Gamma_0(\reals,\reals)$, the claim
  follows.
\end{proof}

Lemma \ref{lem:complexification_and forgetting} implies that, after
inverting $2$, the long exact sequences of Section
\ref{sec:real-versus-complex-1} give rise to the split short exact
sequences we claim to get. 

To relate the K-theory of a complex $C^*$-algebra with the K-theory of
the same $C^*$-algebra, considered as a real $C^*$-algebra, we already
used the following results.

\begin{proposition}\label{prop:compare_real_and_complex}
  Let $\Gamma$ be a discrete group. If $A$ is a $\sigma$-unital
  complex $\Gamma$-$C^*$-algebra (which can also be
  considered as a real $C^*$-algebra) and $B$ is a separable real
  $\Gamma$-$C^*$-algebra, then the inclusion $B\into B_\complexs$ induces a natural
  isomorphisms
  \begin{equation}\label{eq:reals_and_complexs_iso}
   b\colon KK^\Gamma_n(B_\complexs, A) \xrightarrow{\iso}     KKO^\Gamma_n(B, A)
 \end{equation}
\end{proposition}
\begin{proof}
    The isomorphism of \eqref{eq:reals_and_complexs_iso} is given by the
  fact that there is a one to one correspondence already on the level
  of Kasparov triples: since $A$ is $\sigma$-unital, every Hilbert
  $A_\complexs$-module $E$ is a complex vector space, and
  therefore the same is true for the set of bounded operators on
  $E$. Therefore, the real linear maps $B\to\boundedops(E)$ are in
  one-to-one correspondence with the complex linear maps
  $B_\complexs\to\boundedops(E)$. All the other conditions on equivariant
  Kasparov triples, and the equivalence relations are preserved by
  this correspondence. 
All this follows directly by inspecting
  Definitions 2.1 to 2.3 in \cite{MR88j:58123}.
\end{proof}

\sh{Proof of Theorem \ref{theo:general_KK sequence}}

  Special cases of Theorem \ref{theo:general_KK sequence} are well
  known, compare e.g.~\cite[(3.4)]{A1}. We are going to use these
  known results below.

  Following the notation of \cite[Section 2]{A1} and \cite[Section
  7]{Karoubi2}, let $\reals^{1,0}$ be the real line with the
  involution $x\mapsto -x$, and $D^{1,0}$, $S^{1,0}$ the unit disc and
  sphere, respectively, with the induced involution.

  Given any real $C^*$-algebra $A$ and a locally compact space $X$
  with involution $x\mapsto \overline{x}$, following \cite[Section
  6]{Karoubi2} we  define 
  \begin{equation*}
  A(X) := \{ f\colon X\to A\tensor \complexs \mid
  f(\overline{x})=\overline{f(x)}; f(x)\xrightarrow{x\to\infty} 0\}.
\end{equation*}
This is again a real $C^*$-algebra. 

We have the short exact sequence
\begin{equation}\label{eq:small_short_seq}
  0\to \reals(\reals^{1,0})\to \reals(D^{1,0})\to \reals(S^{1,0})\to 0,
\end{equation}
where we identify $\reals^{1,0}$ with the open unit interval in
$D^{1,0}$. This short exact sequence admits a completely positive
cross section, using a cutoff function $\rho\colon D^1\to [0,1]$ with
value $1$ at the boundary and $0$ at the origin to extend functions on
$S^{1,0}$ to the disc.

For an arbitrary real $\Gamma$-$C^*$-algebra $A$, tensoring
\eqref{eq:small_short_seq} with $A$, we get a
short exact sequence (using the canonical isomorphism $A\tensor_\reals
\reals(X)\iso A(X)$)
\begin{equation}\label{eq:big_short_exact}
  0 \to A(\reals^{1,0}) \to A(D^{1,0}) \to A(S^{1,0}) \to 0,
\end{equation}
which again admits a completely positive cross section induced from
the completely positive cross section of \eqref{eq:small_short_seq}. 

Moreover, evaluation at $1$ gives a natural and canonical
$C^*$-algebra \emph{isomorphism}
\begin{equation*}
 \phi\colon A(S^{1,0})\xrightarrow{\iso} A_\complexs.
\end{equation*}
We also note that the evaluation map
\begin{equation*}
   A(D^{1,0}) \to A; f\mapsto f(0)
 \end{equation*}
 is a homotopy equivalence in the sense of KK-theory. The homotopy inverse maps
 $x\in A$ to the constant map with value $x$. In particular, we have
 natural isomorphisms
 \begin{equation*}
   \psi\colon KKO_*^\Gamma(B,A)\iso KKO_*^\Gamma(B,A(D^{1,0})),
 \end{equation*}
 where the maps in both directions are small homomorphisms in the
 sense of Definition \ref{def:small_homomos}. 

  Concerning $A(\reals^{1,0})$, 
by \cite[Paragraph 5, Theorem 7]{MR81m:58075}
  \begin{equation*}
\alpha\colon KKO^\Gamma_n(B,A(\reals^{1,0}))\xrightarrow{\iso} KKO^\Gamma_{n-1}(B,A),
\end{equation*}
where the map $\alpha$ and its inverse are given by exterior
Kasparov product with small elements.

(We remark that the corresponding, but slightly different statements
in \cite[2.5.1]{Schroe} are partly \emph{wrong}, since Schr{\"o}der is
disregarding the gradings.)

  Therefore, the short exact sequence \eqref{eq:small_short_seq}
  induces via the short exact sequence \eqref{eq:big_short_exact}
  for separable real $\sigma$-unital $\Gamma$-$C^*$-algebras $A$ and
  $B$ a long exact sequence in equivariant bivariant real
  K-theory, and we get the following commutative diagram:
\begin{equation*}
    \begin{CD}
@VVV @VVV @VVV\\
KKO_{n+1}^\Gamma(B,A(S^{1,0})) @>{\iso}>{\phi_*}>
KKO_{n+1}^\Gamma(B,A_\complexs)  @<{\iso}<{b}<
KK^\Gamma_{n+1}(B_\complexs,A_\complexs)\\
@V{\delta}VV @V{\delta'}VV @V{\delta''}VV\\
KKO_n^\Gamma(B,A(\reals^{1,0}))  @>>{=}>
KKO_{n}^\Gamma(B,A(\reals^{1,0})) @>{\iso}>{\alpha}>
KKO^\Gamma_{n-1}(B,A)\\
@V{i}VV @V{\rho}VV @V{\chi}VV\\
KKO_n^\Gamma(B,A(D^{1,0})) @<{\iso}<{\psi}<  KKO_n^\Gamma(B,A)  @>>{=}>KKO^\Gamma_{n}(B,A)\\
@V{j}VV @V{i_2}VV @V{c}VV\\
 KKO^\Gamma_n(B,A(S^{1,0})) @>{\iso}>{\phi_*}>
 KKO_n^\Gamma(B,A_\complexs)  @<{\iso}<{b}<
 KK^\Gamma_n(B_\complexs,A_\complexs)\\
@VVV @VVV @VVV\\
  \end{CD}
\end{equation*}

In this diagram, all the horizontal maps are isomorphisms which have
been explained above. The vertical maps in the middle and right row
are induced by the maps of the left row and the horizontal
isomorphisms. The horizontal maps are
small homomorphisms where the inverse is also a small
homomorphism. In particular $\phi_*$ is induced from the
isomorphism of $C^*$-algebras $\phi\colon\reals(S^{1,0})\to\complexs$
(and therefore is small).
Here, $b$ is a somewhat exceptional homomorphism: it is the
identification of KK-groups of Proposition
\ref{prop:compare_real_and_complex} which is true on the level of cycles.

The construction of the long exact sequence in KK-theory implies that
$\delta$, $i$ and $j$ are small homomorphisms (compare
\cite[Theorem 2.5.6]{Schroe}), and $i$ and $j$ are induced from the
maps in the short exact sequence \eqref{eq:small_short_seq}.

Since the compositions of small homomorphisms are small, the same
follows for $\delta'$, $\rho$, $i_2$ and $\chi$. Composing the maps,
we see that $i_2$ is induced from the inclusion
$\reals\into\complexs$, and therefore $c$ is the complexification
homomorphism of Definition \ref{def:complexification}. Since
$\delta''$ equals $\alpha\circ \delta'$ (upto the canonical
identification $b$) we can also consider $\delta''$ as a small
homomorphism.

To identify the small homomorphisms $\tilde\chi\in
KK_1(\reals,\reals)=\integers/2$ and $\delta''$, it suffices to study the case
$A=B=\reals$ and $\Gamma=\{1\}$. Then we get
  \begin{equation}
    \label{eq:long_ex_r_C}
\cdots  \to   K^{n+1}(\complexs) \xrightarrow{\delta''}
KO^{n-1}(\reals)\xrightarrow{\chi}
KO^n(\reals)\xrightarrow{c} K^n(\complexs)\to\cdots
\end{equation}
This exact sequence, including the identification of $\chi$ as
multiplication with the generator of
$KKO^1(\reals)=KKO_1(\reals,\reals)=\integers/2$ and of $\delta''$ as
composition of complex Bott periodicity with ``forgetting the complex
structure'' is already 
established in \cite{A1}. Alternatively, a careful analysis of the
constructions in this special case also identifies $\chi$ and
$\delta''$ without much difficulty. Note that without these
computations  we nevertheless identify $\tilde \chi=\chi$, since
it has to be non-zero by the known K-theory of $\reals$ and
$\complexs$, and there is only one non-zero element in
$KK_1(\reals,\reals)$. Since, in this paper, we don't use the explicit
description of $\delta''$, the main results of this paper are
established without using \cite{A1} (or Remark \ref{rem:calculate_delta}). 

Finally, observe that $\chi$ has additive order
$2$. Moreover, if we iterate $\chi$ three times, we take the Kasparov
product with third power of the generator in $KK_1(\reals,\reals)$
which is zero, and therefore $\chi^3=\chi\circ\chi\circ\chi=0$.

\begin{remark}\label{rem:calculate_delta}
  A possible way to calculate $\delta''$ is the following: the same
  arguments which lead to the K-theoretic exact sequence
  \eqref{eq:long_ex_r_C} give rise to a corresponding sequence in
  K-homology. Using that we know $KKO_*(\reals,\reals)$, this implies that
  $KKO_{-2}(\complexs,\reals)\iso\integers$. Moreover, using the fact
  that ``forgetting the complex structure'' gives a unital ring homomorphism
  $KK_*(\complexs,\complexs)\to KKO_*(\complexs,\complexs)$ (using the
  arguments of Proposition \ref{prop:compare_real_and_complex}) we get
  Bott periodicity $KKO_{-2}(\complexs,\reals)\iso
  KKO_0(\complexs,\reals)$, where the map is given by multiplication
  with the complex Bott periodicity element $\beta\in
  KKO_2(\complexs,\complexs)$. In particular, the
  generators of $KKO_{-2}(\complexs,\reals)$ are products of the
  inverse $\beta^{-1}\in KKO_{-2}(\complexs,\complexs)$ of the complex
  Bott periodicity element with the generators
  of $KKO_0(\complexs,\reals)$. Using the K-homology version of
  \eqref{eq:long_ex_r_C} again (where $\integers\iso KKO_0(\complexs,\reals)\to
  KKO_0(\reals,\reals)\iso\integers$ is induced by the inclusion
  $\reals\into\complexs$), the element ``forgetting the complex
  structure'' $f_\reals$ (defined in \ref{def:forgetting complex structure})
  defines a generator of $KKO_0(\complexs,\reals)$ (since this element
  is mapped to $2$; and we conclude from the long exact sequence that
  this also
  happens to a generator). It follows that the Kasparov product of
  $\beta^{-1}$ with $f_\reals$ is an additive
  generator of $KKO_{-2}(\complexs,\reals)$.

  The exact sequence \eqref{eq:long_ex_r_C} implies that $\delta''$
  can not be divisible. Since it is given by some element of
  $\integers\iso KKO_{-2}(\complexs,\reals)$, it has (up to a sign which we don't
  have to determine) to coincide with the generator, i.e.~the product of
  $\beta^{-1}$ and $f_\reals$.
\end{remark}

\sh{Proof of Theorem \ref{theo:real_versus_complex_BC}}

It remains to prove that the Baum-Connes assembly maps are compatible
with the maps in the long exact sequences. To do this, we have to
recall this assembly (or index) map. It is in the real and complex case
given by the same procedure, which we describe for complex K-theory.

$\mu$ is given as composition of two maps. The first of these is
\begin{equation}\label{eq:descent}
  descent\colon KK_n^\Gamma(C_0(\underline{E}\Gamma),\complexs)\to
  KK_n(C_{red}^*( \Gamma;C_0(\underline{E}\Gamma)), C_{red}^*(\Gamma;\complexs)).
\end{equation}
To be precise, we have to apply this map to $\Gamma$-compact subsets of
$\underline{E}\Gamma$ and then pass to the limit. We avoid this to
simplify notation. By \cite[Theorem 20.6.2]{Bla} (and
\cite[Theorem 2.4.13]{Schroe} for the real case) this descent is
compatible with Kasparov products. Since $\Gamma$ acts trivially on
the right hand factor $\complexs$, $C_{red}^*(\Gamma;\complexs)\iso
\complexs\tensor_\complexs C_{red}^*\Gamma\iso C_{red}^*\Gamma$. Note
that the construction of descend for
trivial actions just amounts to the exterior Kasparov product with the
identity on the level of KK. In other words, if $A$ and $B$ have a
trivial $\Gamma$-action, then
\begin{equation*}
  descent\colon KK(A,B)\to KK(A\tensor C_{red}^*\Gamma,B\tensor C_{red}^*\Gamma)
\end{equation*}
is given by exterior tensor product with the identity.

Since descent is compatible with the intersection product, it follows
that descent commutes with exterior Kasparov product with small
elements in the sense of Definition \ref{def:small_homomos} in
\eqref{eq:descent}.

The second map
\begin{equation}
  \label{eq:mishline}
  KK(C_{red}^*(\Gamma;C_0(\underline{E}\Gamma)),C_{red}^*\Gamma) \to
  KK(\complexs,C_{red}^*\Gamma)
\end{equation}
is given by left Kasparov product with a certain element, the ``Mishchenko
line bundle'', in $KK(\complexs, C_0(\underline{E}\Gamma))$. This also
commutes with exterior Kasparov product with small KK-elements.

Now the Baum-Connes assembly map $\mu$ is the composition of the two
homomorphisms just described, and therefore also commutes with
small homomorphisms.

Since the homomorphisms in the long exact sequences are all small, the
Baum-Connes maps are compatible with them, which is the assertion of
Theorem \ref{theo:real_versus_complex_BC}.

\sh{Variations}

  It is clear that all the arguments given in this section apply in
  exactly the same way in the situations
  described in Theorem \ref{theo:variations}, which is therefore also true.

\Addressesr
\small\url{http://www.uni-math.gwdg.de/schick}

\end{document}